\documentclass[graybox]{svmult}

\usepackage{amsmath,amssymb,bbm}
\usepackage{setspace} 

\usepackage{amsfonts}
\usepackage{mathptmx}
\usepackage{mathptmx}       % selects Times Roman as basic font
\usepackage{helvet}         % selects Helvetica as sans-serif font
\usepackage{courier}        % selects Courier as typewriter font
\usepackage{type1cm}        % activate if the above 3 fonts are
\usepackage{makeidx}         % allows index generation
\usepackage{graphicx}        % standard LaTeX graphics tool
\usepackage{multicol}        % used for the two-column index
\usepackage[bottom]{footmisc}% places footnotes at page bottom

\makeindex             % used for the subject index
\begin{document}

\def\R{\mathbb{R}}
\def \1{1 \mkern -6mu 1} 
\def\N{\mathbb{N}}
\def\E{\mathbb{E}}
\def\P{\mathbb{P}}
\def\F{{\bf F}} 
\def\f{\mathcal{F}}
\def\Z{\mathbb{Z}}
\def\R{\mathbb{R}}
\def\C{\mathbb{C}}
\def\Q{\mathbb{Q}}
\def\L{\mathbb{L}}
\def \e{{\rm e}}
\def \d{{\rm d}}
\def \p{{\mathcal P}}

\title*{Increasing processes and the change of variables formula for non-decreasing functions}
\author{Jean Bertoin and Marc Yor}
\institute{Jean Bertoin \at Institut f\"ur Mathematik, 
Universit\"at Z\"urich, 
Winterthurerstrasse 190, 
CH-8057 Z\"urich, Switzerland, \email{jean.bertoin@math.uzh.ch}
\and Marc Yor \at Institut Universitaire de France and Laboratoire de Probabilit\'es et Mod\`eles Al\'eatoires, 
UPMC, 4 Place Jussieu, 75252 Paris cedex 05, France \email{deaproba@proba.jussieu.fr}}
\maketitle

\abstract {Given an increasing process $(A_t)_{t\geq 0}$, we characterize the right-continuous non-decreasing functions $f: \R_+\to \R_+$
that map $A$ to a pure-jump process. As an example of application, we show for instance that functions with bounded variations belong to the domain of the extended generator of any subordinators with no drift and infinite L\'evy measure.}

\keywords{ Pure jump processes, increasing processes, subordinator, extended infinitesimal generator.}

%\classification{60J30, 60J35}

\hskip 4mm

\section{Introduction}
We make first some simple observations about the composition of non-decreasing functions which partly motivate the present work.
Let $a: \R_+\to \R_+$ be a non-decreasing function, which, as usual, we further assume to be right-continuous with $a(0)=0$. 
 It is well-known that
$a$ can be decomposed canonically into the sum of two non-decreasing functions, $a=a^c+a^d$, where
$a^c$  is continuous and $a^d$ purely discontinuous. 
More precisely, the latter is given by 
$$a^d(t)=\sum_{0< s \leq t} \Delta a(s)\,,$$
where $\Delta a(s)= a(s)-a(s-)$ and the sum above only takes into account the instants $s$ of effective jumps, that is such that $ \Delta a(s)>0$.

We next recall the Stieltjes change of variables formula, a particular case of It\={o}'s formula. Let $f: \R\to \R$ be a ${\mathcal C}^1$ function; then there is the identity
\begin{equation}\label{E0}
f\circ a(t)-f\circ a(0) = \int_0^t f'\circ a(s) \d a^c(s) + \sum_{0< s \leq t}\Delta (f\circ a)(s)\,.
\end{equation}
Note that when $f$ is non-decreasing, (\ref{E0}) specifies the canonical decomposition of the non-decreasing function $f\circ a$. 

We assume henceforth that $a$ is purely discontinuous, i.e.  $a^c\equiv 0$; 
the change of variables formula  (\ref{E0}) then reduces to
\begin{equation}\label{E0'}
f\circ a(t)-f\circ a(0) =\sum_{0< s \leq t}\Delta (f\circ a)(s)\,,
\end{equation}
and it is tempting to think that this identity  still holds when $f$ has only finite variations.
However, this is incorrect in full generality  as we shall now see.

We further suppose that $a$ is strictly increasing; in other words, the Stieltjes measure $\d a$ is purely atomic and its support coincides with $\R_+$. The left-inverse $a^{-1}$ of $a$ is given by
$$a^{-1}(x)=\inf\{y>0: a(y)>x\}\,, \qquad x\geq 0\,;$$
one sometimes calls $a^{-1}$ a Devil's staircase due to the fact that it is a continuous, non-decreasing function that remains constant on the neighborhood of Lebesgue almost every $x>0$. In particular the Stieltjes measure $\d a^{-1}$ is continuous and singular with respect to Lebesgue's measure. Then $a^{-1}\circ a(x)=x$ for all $0\leq x < a(\infty)$, as a consequence  $a^{-1}\circ a$ has no jump at all and the formula (\ref{E0'}) fails for $f=a^{-1}$.

Nonetheless, this note is concerned with the validity of the change of variables formula (\ref{E0'}) in the situation where the deterministic function $a$ above is replaced by an increasing random process $(A_t)_{ t\geq 0}$. Our main result specifies when the latter holds true in terms of a notion of left-accessibility of points. 
This enables us to observe that, for a large class of increasing processes $A$ which includes subordinators with no drift and infinite L\'evy measure, the version of (\ref{E0'}) with $A$ replacing $a$ holds a.s. for any function $f$ with finite variations. 
As an example of application, we deduce that functions with bounded variations belong to the domain of the extended generator of such subordinators.

\section{Increasing processes and left-accessibility of points}

Throughout this note, an increasing process $A=(A_t)_{t\geq 0}$ is  a random process 
 with values in $\R_+$, such that with probability one, $A_0=0$ and the sample path $t \mapsto A_t$ is right-continuous and non-decreasing. 
For every $x>0$, we write
$${L}(x)=\sup\{t\geq 0: A_t<x\}$$
for the last passage time of $A$ below the level $x$. Observe that ${L}(x)$ coincides with the left-limit at $x$ of the right-continuous process $A^{-1}(y)=\inf\{t\geq 0: A_t>y\}$. In particular, ${L}$ is non-decreasing and left-continuous, $A_{{L}(x)}\geq x$ and $A_t<x$ for all $t<L(x)$.

 We then say that $x$ is {\it left-accessible} for a sample path of $A(\omega)$ if
$$A_{{L}(x)-}(\omega)=x\,,$$
that is to say that the sample path takes values strictly less than $x$ which can be arbitrarily close to $x$.
We write ${\mathcal L}(A)$ for the random set of left-accessible points for $A$; note that ${\mathcal L}(A)$ is necessarily closed on the left, i.e. the limit of an increasing sequence of left-accessible points is again left-accessible.  We also define
$$\L(A)=\{x>0: \P(x\in{\mathcal L}(A))>0\}\,.$$
In words, $\L(A)$ is the set of $x>0$ which are left-accessible for $A$ with positive probability.

Finally, let $f:\R_+\to \R_+$ be a right-continuous non-decreasing function.
For every Borel set $B\subseteq \R_+$, we write $f(B)$ for the mass assigned to $B$ by the Stieltjes measure  $\d f$. In particular $f(y)-f(x)=f((x,y])$ for $0\leq x \leq y$.

We now claim the main result of this note. 

\begin{theorem} \label{T1}
{\rm (i)} If $f(\L(A))=0$, then we have
$$f(A_t)-f(0)= \sum_{0< s \leq t} \left(f(A_s)-f(A_{s-})\right)$$
 for all $t\geq 0$, a.s.
 
 {\rm (ii)} If $f(\L(A))>0$, then for $t>0$ sufficiently large, the strict inequality
$$f(A_t)-f(0)> \sum_{0< s \leq t} \left(f(A_s)-f(A_{s-})\right)$$
 holds with positive probability.
\end{theorem}

 \proof We start by pointing out that Fubini-Tonelli Theorem yields the identity
 \begin{eqnarray*}
 \E(f({\mathcal L}(A))) &=&\E\left( \int_{(0,\infty)} \d f(x) \1_{{\mathcal L}(A)}(x)\right) \\
& =&  \int_{(0,\infty)} \d f(x) \P(x\in {\mathcal L}(A)) \,.
 \end{eqnarray*}
 Thus $f(\L(A))=0$ if and only if $f({\mathcal L}(A))=0$ a.s.

(i) We assume that $f(\L(A))=0$ and denote by $\Lambda$ the event that $f({\mathcal L}(A))=0$, so $\P(\Lambda)=1$. We pick $\omega\in \Lambda$ and observe that for every $x>0$ which is not
left-accessible for $A(\omega)$, there is the identity
$$\1_{[x,\infty)}(A_t(\omega))= \sum_{0< s \leq t}\left(\1_{[x,\infty)}(A_s(\omega)) - \1_{[x,\infty)}(A_{s-}(\omega)) \right)\,.$$
Indeed the left-hand side equals $1$ if $t\geq {L}(x)$ and $0$ otherwise. Note that since $A(\omega)$ is non-decreasing,  the summand in the right-hand side is $0$ for all $s>0$, except for $s={L}(x)$ where it equals $1$ since $x$ is not left-accessible.
We also point out that the sum in the right-hand side {\it de facto}  only involves instants $s>0$ at which $A(\omega)$ is discontinuous, and since the set of discontinuities of $A(\omega)$ is at most countable, the counting measure on that set is sigma-finite.

We now integrate the preceding equality for all $x>0$ which are not left-accessible for $A(\omega)$, with respect to the Stieltjes measure $\d f$. Recall that $\d f$ gives no mass to the set of left-accessible points; we get that the left-hand side equals $f(A_t(\omega))-f(0)$. We apply the Fubini-Tonelli theorem on the right-hand side, which is legitimate by the observation which was made  above. We obtain
\begin{eqnarray*}
& &\sum_{0< s \leq t}\int_{(0,\infty)} \left(\1_{[x,\infty)}(A_s(\omega)) - \1_{[x,\infty)}(A_{s-}(\omega)) \right)  \d f(x) \\
&=& \sum_{0< s \leq t} \left(f(A_s(\omega))-f(A_{s-}(\omega))\right)\,,
\end{eqnarray*}
which proves our claim.

(ii) We now assume  that $f(\L(A))>0$ and denote by $\Lambda^{\rm c}$ the event that $f({\mathcal L}(A))>0$, so $\P(\Lambda^{\rm c})>0$.
We pick $\omega\in \Lambda^{\rm c}$ and observe that for every $x>0$ which is 
left-accessible for $A(\omega)$ and $t\geq {L}(x)$,
we have
$$\1_{[x,\infty)}(A_t(\omega))= 1 > 0 =\sum_{0< s \leq t}\left(\1_{[x,\infty)}(A_s(\omega)) - \1_{[x,\infty)}(A_{s-}(\omega)) \right)\,.$$
It then follows from arguments based on the Fubini-Tonelli theorem similar to those in (i) that
$$f(A_t(\omega))-f(0)> \sum_{0< s \leq t} \left(f(A_s(\omega))-f(A_{s-}(\omega))\right)$$
whenever $t$ is sufficiently large. 
{\hfill $\Box$}  \endproof

We stress that for a large class of increasing processes $A$, the set $\L(A)$ is empty and therefore the change of variables formula holds a.s. for any  function with finite variations,
that is which can be expressed as the difference of two non-decreasing right-continuous functions. 
Indeed, observe that if $x\in \L(A)$, then the probability that there exists $s>0$ with $A_{s-}=x$ must be positive. Thus $\L(A)=\varnothing$ whenever every single point $x>0$
is {\it polar} for $A$, that is  if 
$$\P(\exists s>0: A_s=x \hbox{ or } A_{s-}=x)=0\,.$$
According to an important result due to Kesten \cite{Ke} and Bretagnolle \cite{Bre} (see also Andrew \cite{An} for a more elementary argument),  polarity of single points holds for any subordinator with no drift and infinite L\'evy measure.  Further examples  can then be constructed from driftless subordinators, e.g. by strictly increasing mapping, change of time or locally equivalent change of probability measures. For instance, the well-know correspondence due to Lamperti \cite{La} connecting positive self-similar Markov processes and L\'evy processes, entails that single points are polar for any strictly increasing self-similar Markov process which is further pure-jump.

\section{Some applications}
In this section, we  present some direct  consequences of Theorem \ref{T1}.

First recall the following fact that has been pointed out by Kingman \cite{Ki} in the setting of subordinators. 
A non-decreasing function $a: \R_+\to \R_+$ is purely discontinuous if and only if its closed range,
$${\mathcal R}(a)=\{x\geq 0: x=a(s) \hbox{ or } x=a(s-) \hbox{ for some }s\geq 0\}\,,$$
has zero Lebesgue's measure, i.e. $|{\mathcal R}(a)|=0$. 
Indeed, for every $t\geq 0$, the complementary set $[a(0),a(t)]\backslash {\mathcal R}(a)$ is  open and  can be expressed as the union of the disjoint open intervals $(a(s-),a(s))$ for $0<s\leq t$ (of course, such an interval is non-empty if and only if $s$ is a  discontinuity point of $a$, which forms a set that is at most countable). 
When $|{\mathcal R}(a)|=0$, this open set has Lebesgue's measure $a(t)-a(0)$, which yields the identity
$$a(t) -a(0) = \sum_{0<s \leq t} (a(s)-a(s-))\,.$$ 
Conversely, when $a$ is purely discontinuous, the identity above holds for all $t\geq 0$ and this entails that $|{\mathcal R}(a)|=0$. 

Combining with Theorem \ref{T1}, we immediately get the following. 

\begin{corollary}Ê\label{C1} Let $f:\R+\to \R_+$ be a right-continuous non-decreasing function.

\noindent {\rm (i)}
 If $f(\L(A))=0$, then $|{\mathcal R}(f(A))|=0$ a.s. 
 
\noindent  {\rm (ii)}
Conversely, if $f(\L(A))>0$ and we assume further that $f$ is continuous, then $\P(|{\mathcal R}(f(A))|>0)>0$. 
\end{corollary}
{\bf Remark.} The continuity of $f$ is required in the second part of the statement to ensure that the only jump times of the process $f(A)$ are those of $A$.

We next turn our attention to stochastic calculus. 
Consider some filtered probability space $(\Omega, \f, (\f_t), \P)$ which fulfills the usual conditions. Let $(X_t)_{t\geq 0}$ be a c\`adl\`ag semi-martingale and $(A_t)_{t\geq 0}$ a c\`adl\`ag increasing (adapted) process which is purely discontinuous. Consider further a function $h: \R_+\to \R$ of class ${\mathcal C}^1$. Then it is well-known that $h(A)$ is a process with finite variations, and the following integration by parts formula holds:
\begin{eqnarray}\label{E1}
& &X_t h(A_t) - X_0h(A_0)\nonumber \\
&=& \int_0^t h(A_{s-}) \d X_s + \sum_{0< s \leq t} X_{s-}\left(h(A_s)-h(A_{s-})\right) +  \sum_{0< s \leq t} \Delta X_s\left(h(A_s)-h(A_{s-})\right) \nonumber \\
&=& \int_0^t h(A_{s-}) \d X_s + \sum_{0< s \leq t} X_s\left(h(A_s)-h(A_{s-})\right).
\end{eqnarray}

We now consider  a c\`adl\`ag  function $k: \R_+\to \R$ with finite variations, that is that $k$ can be expressed as the difference of two right-continuous non-decreasing functions.
For every Borel set $B\subset \R_+$, we write $|k|(B)$ for the mass given to $B$ by the total-variation measure $|\d k|$. If we assume that $|k|(\L(A))=0$, then it follows from Theorem \ref{T1} that $k(A)$ is a purely discontinuous process with finite variations, and the classical integration by parts formula yields 
$$X_t k(A_t) -X_0k(0)= \int_0^t k(A_{s-}) \d X_s + \sum_{0< s \leq t} X_{s}\left(k(A_s)-k(A_{s-})\right) \,.$$
In particular, when $\L(A)=\varnothing$, and {\it a fortiori}  when single points are polar for $A$, the integration by parts formula (\ref{E1}) holds  when one merely assumes that the function $h$ has finite variations.

We now conclude this note with an application to the infinitesimal generator of a subordinator. We consider 
a subordinator  $(S_t)_{t\geq 0}$  with no drift and infinite L\'evy measure. 
That is, $S$ is a random process with values in $\R_+$, with independent and stationary increments, and we assume that  its sample path is right-continuous, strictly increasing, and has no continuous component. According to the L\'evy-It\={o} decomposition, the Stieltjes measure $\d S$ is expressed in the form
$$\d S = \sum_{t\geq 0} \Delta S_t \delta_t$$
where $ \Delta S_t=S_t-S_{t-}$ denotes the possible jump at time $t$.
More precisely, the jump process $\Delta S$ is a Poisson point process whose intensity is known as the L\'evy measure of $S$.
We refer to \cite{BeSF} for background. 

Recall that $S$ is a Feller process on $\R_+$; we write $G: {\mathcal D}\to {\mathcal C}_0$ for its infinitesimal generator, where ${\mathcal D}$ is its domain and ${\mathcal C}_0$ the space of continuous functions on $\R_+$ with limit $0$ at $+\infty$. It is well-known that if ${\mathcal C}^1_0$ denotes the subspace of functions $g\in {\mathcal C}_0$ which are continuously differentiable with $g'\in {\mathcal C}_0$, then $ {\mathcal C}^1_0\subset {\mathcal D}$ and there is the identity
$$Gg(x)=\int_{(0,\infty)} \left( g(x+y)-g(x)\right) \Pi(\d y)\,, \qquad x\geq 0\,,$$
where $\Pi$ denotes the L\'evy measure of $S$. Recall also that for every $g\in{\mathcal D}$ and $x\geq 0$, the process
$$g(x+S_t) - \int_0^t Gg(x+S_s) \d s $$
is a martingale, and conversely, this martingale property together with the assumption that both $g$ and $Gg$ are in ${\mathcal C}_0$, characterize the infinitesimal generator. Further, a measurable function $h: \R_+\to \R$ belongs to the  domain ${\mathcal D}_e$ of the {\it extended} generator if there exists a measurable function $b: \R_+\to \R$  such that for every $x\geq 0$ and $t\geq 0$
$$\int_0^t |b(x+S_s)| \d s < \infty \qquad \hbox{a.s.}$$
and 
$$h(x+S_t) - \int_0^t b(x+S_s) \d s $$
is a martingale. By slightly abusing notation (because $b$ is not unique), we then write $G_eh=b$.

\begin{corollary} Let $f: \R_+\to \R$ be a non-decreasing right-continuous function, which we further assume to be bounded. Then $f\in {\mathcal D}_e$ and 
$$G_ef(x)=\int_{(0,\infty)} \left( f(x+y)-f(x)\right) \Pi(\d y)\,, \qquad x\geq 0\,,$$
where $\Pi$ denotes the L\'evy measure of $S$.
\end{corollary}
Note that $G_ef(x)$ is well-defined in $[0,\infty]$ (and thus possibly infinite) since $f$ is non-decreasing.
\proof   Recall that single points are polar for $S$ according to the result of Kesten \cite{Ke}, and thus Theorem \ref{T1}(i) applies to any non-decreasing function $f$. 
Replacing the function $f$ there by the function $f(x+\cdot)$, we get
\begin{eqnarray*}
f(x+S_t)-f(x) &=& \sum_{0< s \leq t} \left(f(x+S_s)-f(x+S_{s-})\right) \\
&=& \sum_{0< s \leq t} \left(f(x+S_{s-}+\Delta S_s)-f(x+S_{s-})\right) \,.
\end{eqnarray*}
Recall that the jump process $\Delta S $ of $S$ is a Poisson point process with intensity $\Pi$. It follows that the predictable compensator of the increasing process $f(x+S_t)-f(x)$ is
$$\int_0^t\d s \int_{(0,\infty)}\Pi(\d y) \left( f(x+S_s+y)-f(x+S_s)\right)=\int_0^t G_ef(x+S_s) \d s .$$
In other words, 
$$f(x+S_t)-f(x)- \int_0^t G_ef(x+S_s) \d s$$
is a martingale. {\hfill $\Box$}  \endproof

Of course, by linearity, we deduce that functions $k: \R_+\to \R$ with bounded variations, i.e. which can be expressed as the difference of two bounded non-decreasing functions, are also in the domain of the extended generator, and then 
$$G_ek(x)=\int_{(0,\infty)} \left( k(x+y)-k(x)\right) \Pi(\d y)\,,$$
where we agree for instance that $G_ek(x)=0$ when the integral above is not well-defined.

\end{document}